\documentclass[12pt,twoside,a4paper]{article}
\usepackage{amsmath,amssymb}
\usepackage{fullpage}
\usepackage{pslatex}
\usepackage{url}

\newcommand{\Eulerphi}{\varphi}
\newcommand{\dx}{\mathrm{d}}

\newtheorem{Lemma}{Lemma}
\newenvironment{Proof}[1][Proof]{\par\noindent\textbf{#1.}~}
  {\hfill$\square$\smallskip\par}

\begin{document}

\title{On the constant in the Mertens product \\
       for arithmetic progressions. \\
       II. Numerical values}
\author{A.~LANGUASCO and A.~ZACCAGNINI}

\date{December 11, 2007}

\maketitle

\begin{abstract}
We give explicit numerical values with 100 decimal digits for the
constant in the Mertens product over primes in the arithmetic
progressions $a \bmod q$, for $q \in \{3$, \dots, $100\}$ and $(a, q) = 1$.
\end{abstract}

\section{Introduction}

In our recent paper \cite{LanguascoZaccagnini2007} we found a new
expression for the constant $C(q, a)$ defined implicitly by
\begin{equation}
\label{def-C}
  P(x; q, a)
  =
  \prod_{\substack{p \le x \\ p \equiv a \bmod q}}
    \Bigl( 1 - \frac1p \Bigr)
  =
  \frac{C(q, a)}{(\log x)^{1 / \Eulerphi(q)}}
  (1 + o(1))
\end{equation}
as $x \to +\infty$, where, here and throughout the present paper,
$q \ge 3$ and $a$ are fixed integers with $(a, q) = 1$, and $p$
denotes a prime number.
When $q \in \{1$, $2\}$ the value of $C(q, a)$ can be deduced from the
classical Mertens Theorem.
In particular, we proved that
\begin{equation}
\label{value-C}
  C(q, a)^{\Eulerphi(q)}
  =
  e^{-\gamma}
  \prod_p
    \Bigl( 1 - \frac1p \Bigr)^{\alpha(p; q, a)}
\end{equation}
where $\alpha(p; q, a) = \Eulerphi(q) - 1$ if $p \equiv a \bmod q$ and
$\alpha(p; q, a) = -1$ otherwise, and $\gamma$ is the Euler constant.
The infinite product is convergent, though not absolutely, by the
Prime Number Theorem for Arithmetic Progressions.

In our paper \cite{LanguascoZaccagnini2007d} we proved that the
constants $C(q, a)$ satisfy some interesting identities but,
unfortunately, these are not suitable for numerical computations.
Here we derive further identities, involving Dirichlet $L$-functions,
that enable us to compute numerically the values of $C(q, a)$ with
many digits for comparatively small $q$.
Details of these identities are given in \S\ref{sec:numer-comp}, and
the results of our numerical computations are collected in
\S\ref{sec:numer-values}: some sample values, truncated to $40$
decimal digits, are shown in Tables~\ref{firsttable}--\ref{thirdtable}.
Finch \cite{Finch2007} has done some numerical work in the case
$q \in \{ 3$, $4\}$.

We would like to express our warmest thanks to Henry Cohen for an
illuminating discussion.

\section{Theoretical framework}
\label{sec:numer-comp}

In this section we concentrate on the numerical computation of the
values of the constant $C(q, a)$ for comparatively small values of
$q$, starting from our formula \eqref{value-C}, and give the
theoretical framework for the results in \S\ref{sec:numer-values}.
We adhere to the notations in the books by Henri Cohen
\cite{Cohen2007a,Cohen2007b}.

We will use the following convention: for any real positive constant
$A$ and for any Dirichlet $L$ function, we write
\[
  L_A(\chi, s)
  =
  \prod_{p > A} \Bigl( 1 - \frac{\chi(p)}{p^s} \Bigr)^{-1},
\]
and do similarly for other Euler products.
We want to compute
\[
  \Eulerphi(q)
  \log C(q, a)
  =
  -\gamma
  + \log \frac{q}{\Eulerphi(q)}
  -
  \sum_{\substack{\chi \bmod q \\ \chi \ne \chi_0}}
    \overline{\chi}(a)
    \sum_{m \ge 1} \frac1m \sum_p \frac{\chi(p)}{p^m}
\]
Notice that the last sum over $p$ is $\sim \chi(2) 2^{-m}$ when $m$ is
large.
We compute the sum over $p$ by M\"obius inversion;
let $A$ be a fixed positive constant: then 
\[
  \sum_p \frac{\chi(p)}{p^m}
  =
  \sum_{p \le A q} \frac{\chi(p)}{p^m}
  +
  \sum_{k \ge 1}
    \frac{\mu(k)}k \log(L_{A q}(\chi^k, km)).
\]
Therefore
\begin{equation}
\label{fundamental}
  \Eulerphi(q)
  \log C(q, a)
  =
  -\gamma
  +
  \log \prod_{p \le A q}
    \Bigl( 1 - \frac1p \Bigr)^{\alpha(p; q, a)}
  -
  \sum_{\substack{\chi \bmod q \\ \chi \ne \chi_0}}
    \overline{\chi}(a)
    \sum_{m \ge 1}
      \frac1m
      \sum_{k \ge 1}
        \frac{\mu(k)}k \log \bigl( L_{A q} (\chi^k, k m) \bigr).
\end{equation}
Grouping the terms with the same value of $k m$, we see that the last
part is
\[
  \sum_{m \ge 1}
    \frac1m
    \sum_{k \ge 1}
      \frac{\mu(k)}k \log \bigl( L_{A q} (\chi^k, k m) \bigr)
  =
  \sum_{n \ge 1}
    \frac1n
    \sum_{k \mid n} \mu(k) \log \bigl( L_{A q} (\chi^k, n) \bigr)
\]
Notice that the Riemann zeta function is never computed at $s = 1$ in
\eqref{fundamental}, since $k m = 1$ implies $k = 1$, and this in its
turn implies $\chi^k = \chi = \chi_0$.
For $n > 1$ we use
\begin{equation}
\label{errortail}
  \bigl| \log \bigl( L_{A q} (\chi^k, n) \bigr) \bigr|
  \le
  \frac1{(n - 1) (A q)^{n - 1}}.
\end{equation}
This inequality is a consequence of the following Lemma.
We remark that a stronger result is valid for small $n$, but the
simple bound below suffices for our application.

\begin{Lemma}
Let $\chi \bmod q$ be any character and $n \ge 2$ be an integer.
If $B \ge 1$ is an integer then
\[
  \bigl| \log \bigl( L_B(\chi, n) \bigr) \bigr|
  \le
  \frac{B^{1 - n}}{n - 1}.
\]
\end{Lemma}

\begin{Proof}
By the triangle inequality
\[
  \bigl| \log \bigl( L_B(\chi, n) \bigr) \bigr|
  =
  \Bigl|
    \sum_{p > B}
      \sum_{m \ge 1} \frac{\chi^m(p)}{m p^{mn}}
  \Bigr|
  \le
  \sum_{p > B}
    \sum_{m \ge 1} \frac1{m p^{mn}}
  \le
  \sum_{k > B} \frac1{k^n}
  \le
  \int_B^{+\infty} \frac{\dx t}{t^n}
  =
  \frac{B^{1 - n}}{n - 1},
\]
as required.
\end{Proof}

\noindent
We have thus reduced the task of the computation of $\log(C(q, a))$ to
computing $\log \bigl( L_{A q} (\chi^k, n) \bigr)$ to $100$ decimal
places, say.
In what follows we denote by $\chi$ a generic Dirichlet character
$\bmod q$ and by $n \ge 1$ an integer.

\paragraph{First step}
We write
\[
  L_{A q}(\chi, n)
  =
  L(\chi, n)
  \prod_{p \le A q} \Bigl( 1 - \frac{\chi(p)}{p^n} \Bigr)
\]
for a convenient value of $A$.

\paragraph{Second step: reduction to primitive characters}
Assume now that $\chi \bmod q$ is induced by $\chi_f \bmod f$, where $f$
is the conductor of $\chi$.
Then we have the identity
\[
  L(\chi, n)
  =
  L(\chi_f, n)
  \prod_{p \mid q}
    \Bigl( 1 - \frac{\chi_f(p)}{p^n} \Bigr).
\]
In particular, we recall that if $\chi = \chi_0 \bmod q$ then
\[
  L(\chi_0, n)
  =
  \zeta(n)
  \prod_{p \mid q}
    \Bigl( 1 - \frac{1}{p^n} \Bigr).
\]

\paragraph{Third step: first case}
Now assume that $\chi$ is a primitive character modulo $f$ and that
$\chi(-1) = (-1)^n$.
Then, by Proposition~10.2.4 of Cohen \cite{Cohen2007b} we have the
explicit formula
\begin{equation}
\label{expl-formula-L}
  L(\chi, n)
  =
  \frac12 (-1)^{n - 1 + (n + e) / 2}
  W(\chi) \sqrt{f}
  \Bigl( \frac{2 \pi}f \Bigr)^n
  \frac{\overline{B_n(\chi)}}{n!},
\end{equation}
where $W(\chi)$ denotes the \emph{root number} of $\chi$ (see
Definition~2.2.25 in \cite{Cohen2007a}), $e = 0$ if $\chi$ is even and
$e = 1$ if $\chi$ is odd, and $B_n(\chi)$ denotes the
$\chi$-\emph{Bernoulli number} which, in its turn, is defined by means
of the $n$-th \emph{Bernoulli polynomial} $B_n(x)$ (see
\cite{Cohen2007b}, Definition~9.1.1), as follows
\[
  B_n(\chi)
  =
  f^{n - 1}
  \sum_{a = 0}^{f - 1} \chi(a) B_n \Bigl( \frac af \Bigr).
\]
This definition is valid both for primitive and imprimitive
characters.
This is the last identity of Proposition~9.4.5 in Cohen \cite{Cohen2007b}.

\paragraph{Third step: second case}
If $\chi$ is non-principal and $\chi(-1) = (-1)^{n + 1}$, there are two
possibilities.

\begin{itemize}

\item
Use the $\chi$-Euler-MacLaurin summation formula (the number of steps
is proportional to $q$, but all terms are elementary); see Cohen
\cite{Cohen2007b}, Corollary~9.4.18;

\item
Use the functional equation, which is valid if $\chi$ is primitive:
this would take a smaller number of steps, of the order
$\asymp \sqrt{q} \log q$, but it needs the computation of the
incomplete $\Gamma$ function.

\end{itemize}

For $q$ small, we use the Euler-MacLaurin summation formula.
When computing $L(\chi, n)$ with $n$ large, the functional equation
does not take into account the fact that
$L(\chi, n) = 1 + \chi(2) 2^{-n} +$ very much smaller terms.

When using the Euler-MacLaurin formula we take a multiple $N$ of $q$
and for $\Re(s) > 1$ write
\begin{equation}
\label{L-decomp}
  L(\chi, s)
  =
  \sum_{r < N} \frac{\chi(r)}{r^s}
  +
  B_0(\chi)
  \frac{N^{1 - s}}{s - 1}
  -
  \frac1{N^s}
  \sum_{j = 1}^T
    \frac{(-1)^{j - 1} B_j(\chi)}{j!}
    \frac{s (s + 1) \cdots (s + j - 2)}{N^{j - 1}}
  +
  R(T),
\end{equation}
where
\begin{align*}
  R(T)
  &=
  -\frac1{T!}
  s (s + 1) \cdots (s + T - 1)
  \int_N^{+\infty}
    B_T(\chi^-, \{ t \}_{\chi}) \frac{\dx t}{t^{s + T}}, \\
  B_T(\chi^-, \{ t \}_{\chi})
  &=
  f^{T-1}
  \sum_{r \bmod f}
    \chi^-(r) B_T \Bigl( \Bigl\{ \frac{t + r}f \Bigr\} \Bigr)
\end{align*}
and $\chi^-(n) = \chi(-n)$; see the Definitions~9.4.2 and~9.4.10
in \cite{Cohen2007b}.
The asymptotic series above is \emph{not} convergent: we take terms
until $R(T)$ reaches a small minimum, before it starts growing again.

Notice that $B_0(\chi) = 0$ in \eqref{L-decomp} for non-principal
$\chi$ by Proposition 9.4.5 of Cohen \cite{Cohen2007b} and the remarks
following it immediately.
This is indeed crucial for the rapidity of convergence.

When $\chi(-1) = (-1)^n$ we use \eqref{expl-formula-L} to estimate
$B_n(\chi) \asymp (q / 2 \pi)^n$.
If $\chi(-1) = (-1)^{n + 1}$ then $B_n(\chi) = 0$.

\paragraph{Computation of the root number}
If $\chi$ is a primitive character modulo $q$, then the \emph{root
number} $W(\chi)$ is defined by means of
\[
  W(\chi)
  =
  \frac{\tau(\chi)}{\sqrt{q} i^e}
  \qquad\text{where}\qquad
  \chi(-1) = (-1)^e
  \quad\text{and $e \in \{0$, $1\}$},
\]
and $\tau(\chi) = \sum_{r=1}^q \chi(r) e(r/q)$ is the Gauss sum.
It is well known that $|W(\chi)| = 1$.
If $\chi^2 = \chi_0$ then $\chi$ is a Legendre symbol and
$W(\chi) = 1$.

For $q$ small, this is alright.
For $q$ large, we use the functional equation which is valid for
primitive $\chi$, introduce
\[
  c(\chi)
  =
  \sum_{n \ge 1} \chi(n) e^{-\pi n^2}
\]
and notice that
\[
  W(\chi)
  =
  \frac{c(\chi)}{i^e \overline{c(\chi)}}.
\]

\section{Description of the computer program}
\label{sec:numer-values}

First of all we need to generate the complete set of Dirichlet
characters $\bmod q$ and also to compute their orders and conductors
and whether they are primitive or not.
To this end we follow the argument in \S4 of Davenport
\cite{Davenport2000}: we first generate the characters for any
$p^\alpha \mid q$, paying particular attention to the case when $q$ is
an even integer, and then we build by multiplication the characters to
the modulus $p_1^\alpha p_2^\beta$ with $p_1 \ne p_2$ and $p_1$,
$p_2 \mid q$.
To compute the order and the primitivity of this character we use
Proposition~2.1.34 of \cite{Cohen2007a}.
The conductor of a character is obtained using the necessary and
sufficient condition described in Lemma~2.1.32 of \cite{Cohen2007a}.

In order to evaluate \eqref{fundamental} using a computer program we
have to truncate the sums over $k$ and $m$ and to estimate the error
we are introducing.
Let $M$, $K > 1$ be two integers.
We have
\allowdisplaybreaks
\begin{align*}
  \log \prod_{p > A q}
    \Bigl( 1 - \frac1p \Bigr)^{\alpha(p; q, a)}
  &=
  -
  \sum_{\substack{\chi \bmod q \\ \chi \ne \chi_0}}
    \overline{\chi}(a)
    \sum_{1\leq m \leq M} \frac1m \sum_{p>Aq} \frac{\chi(p)}{p^m}
  -
  \sum_{\substack{\chi \bmod q \\ \chi \ne \chi_0}}
    \overline{\chi}(a)
    \sum_{m > M} \frac1m \sum_{p>Aq} \frac{\chi(p)}{p^m} \\
  &=
  -
  \sum_{\substack{\chi \bmod q \\ \chi \ne \chi_0}}
    \overline{\chi}(a)
    \sum_{1\leq m \leq M} \frac1m
    \sum_{1\leq k \leq K}
    \frac{\mu(k)}k \log(L_{A q}(\chi^k, km)) \\
  &-
  \sum_{\substack{\chi \bmod q \\ \chi \ne \chi_0}}
    \overline{\chi}(a)
    \sum_{1\leq m \leq M} \frac1m
    \sum_{k> K}
    \frac{\mu(k)}k \log(L_{A q}(\chi^k, km))\\
  &-
  \sum_{\substack{\chi \bmod q \\ \chi \ne \chi_0}}
    \overline{\chi}(a)
    \sum_{m > M} \frac1m \sum_{p>Aq} \frac{\chi(p)}{p^m} \\
  &=
  - S(q,a) - E_1(q,a,A,K) - E_2(q,a,A,M),
\end{align*}
say.
Using \eqref{errortail} and the trivial bound for $\chi$, it is easy
to see that
\[
  \left\vert E_1(q,a,A,K) \right \vert
  \leq
  \frac{2 A q (\Eulerphi(q) - 1)}{2 K (A q - 1) \Bigl[(Aq)^{K}-1\Bigr]}
\]
and
\[
  \left\vert E_2(q,a,A,M) \right\vert
  \leq
  \frac{A q (\Eulerphi(q) - 1)}{M (M - 1) (A q - 1) (Aq)^{M}}.
\]
In order to ensure that $S(q, a)$ is a good approximation of $C(q, a)$
it is sufficient that $A q$, $K$ and $M$ are sufficiently large.
Setting $A q = 9600$ and $K = M = 26$ yields the desired $100$ correct
decimal digits.

Now we have to consider the error we are introducing during the
evaluation of the Dirichlet $L$-functions that appear in $S(q, a)$.
Notice that in the case involving the Bernoulli numbers we use an
exact formula: hence we just need to evaluate the error introduced by
the $R(T)$ term in the Euler-McLaurin summation formula
\eqref{L-decomp}.
In fact the Euler-McLaurin summation formula is used in about $1/4$ of
the total cases but we are now just looking for an upper bound and so
we will sum $\vert R(T) \vert $ over $m \leq M$ and $k \leq K$.

Assume now that $T\geq 2$ is an even integer and $q \mid N$.
For any non-principal character $\chi^k \bmod{q}$ \eqref{L-decomp}
implies that
\[
  L_{T,N}(\chi^k,km)
  =
  \sum_{r < N} \frac{\chi^k(r)}{r^{km}}
  -
  \frac1{N^{km}}
  \sum_{j = 1}^T
    \frac{(-1)^{j - 1} B_j(\chi^k)}{j!}
      \frac{km (km + 1) \cdots (km + j - 2)}{N^{j - 1}}
\]
and hence we get
\[
  L_{A q}(\chi^k, km)
  =
  \Pi \left(L_{T,N}(\chi^k,km) - E_3(q,m,k,N,T, \chi^k) \right),
\]
where $\Pi$ denotes the finite products we wrote in the first and
second step of \S\ref{sec:numer-comp}.
Moreover it is clear that
\[
  \left\vert E_3(q,m,k,N,T,\chi^k) \right\vert
  \le
  \frac{km(km+1)\dotsm(km+T-1)}{T!}
  \int_N^{+\infty}
    \vert B_T((\chi^k)^{-}, \{t\}_{\chi^k}) \vert \ t^{-km-T} \dx t.
\]
Hence
\[
  \left\vert
  \log(L_{A q}(\chi^k, km))
  -
  \log \left(\Pi \cdot L_{T,N}(\chi^k,km)\right)
  \right\vert
  \leq\left\vert
  \frac{ E_3(q,m,k,N,T,\chi^k)}{L_{T,N}(\chi^k,km)}
  \right\vert
\]
and the total error arising in the computation of the Dirichlet
$L$-functions can be obtained summing the previous estimate over $m$
and $k$.
For $T$ even, trivial estimates and Proposition 9.1.3 of
\cite{Cohen2007b} imply that
\begin{align*}
  \left\vert B_T((\chi^k)^{-}, \{t\}_{\chi^k}) \right \vert
  &\le
  f^{T-1}
  \sum_{r \bmod f}
    \left\vert
      \chi^k(-r) B_T \left(\left\{\frac{t+r}{f} \right\}\right)
    \right\vert \\
  &\le
  f^{T-1}
  \sum_{r \bmod f}
    \left\vert
      \sum_{j=0}^T
        \binom{T}{j} B_j
        \left\{\frac{t+r}{f} \right\}^{T-j}
    \right\vert
  \le
  f^T B_T,
\end{align*}
where $f \mid q$ is the conductor of $\chi^k$ and $B_T$ is the $T$-th
Bernoulli number.
Hence we obtain
\begin{align*}
  \left\vert E_3(q,m,k,N,T,\chi^k) \right\vert
  &\le
  \frac{km(km+1)\dotsm(km+T-1)q^T B_T}{T!}
  \frac{N^{1-km-T}}{km+T-1} \\
  &=
  \frac{q^T B_T}{T!} (km)\dotsm(km+T-2) N^{1-km-T}.
\end{align*}
Moreover, let
\[
  U(q, M, K, N, T)
  =
  \min_{\substack{\chi \bmod q \\ \chi \neq \chi_0}}
  \min_{\substack{1\leq k \leq K \\ 1\leq m \leq M}}
  \vert L_{T,N}(\chi^k,km) \vert.
\]
The total error arising in the computation of the Dirichlet
$L$-functions is therefore
\begin{align*}
  \left\vert E_4(q,a,M,K,N,T) \right\vert
  &\leq
  \frac{(\Eulerphi(q)-1)q^T B_T}{U(q,M,K,N,T)}
  \sum_{1\leq m \leq M} \frac1m
     \sum_{1\leq k \leq K}
       \frac{1}{k}
       \frac{(km)\dotsm(km+T-2)}{T!} N^{1-km-T} \\
  &=
  \frac{(\Eulerphi(q)-1)q^T B_T}{ U(q,M,K,N,T)T!}
  \sum_{1\leq m \leq M}
     \sum_{1\leq k \leq K} (km+1)\dotsm(km+T-2) N^{1-km-T} \\
  &\leq
  \frac{(\Eulerphi(q)-1)(KM+T-2)^{T-2}q^T B_T}{ U(q,M,K,N,T)N^{T-1}T!}
  \sum_{1\leq m \leq M}
    \sum_{1\leq k \leq K} N^{-km} \\
  &\leq
  \frac{2(\Eulerphi(q)-1)(KM+T-2)^{T-2}q^T B_T}{(N-1) U(q,M,K,N,T)N^{T-1}T!}.
\end{align*}
Letting
\[
  \widetilde{C}(q,a)
  =
  \left(
    e^{-\gamma}
    \prod_{p \le A q}
      \Bigl( 1 - \frac1p \Bigr)^{\alpha(p; q, a)}
    \exp(-S(q,a))
  \right)^{1/\Eulerphi(q)}
\]
and collecting the previous estimates, we have that
\begin{align*}
  \Bigl\vert
    C(q,a)
    -
    \widetilde{C}(q,a)
  \Bigr\vert
  &\leq
  \widetilde{C}(q,a)
  \left\vert
    \exp
    \left(
    -\frac{E(q,a,A,M,K,N,T)}{\Eulerphi(q)}
    \right)
    -1
  \right\vert \\
  &\le
  \widetilde{C}(q,a)
  \frac{|E(q,a,A,M,K,N,T)|}{\Eulerphi(q)},
\end{align*}
where $E(q,a,A,M,K,N,T)$ denotes
$E_1(q,a,A,K) + E_2(q,a,A,M) + E_4(q,a,M,K,N,T)$.

Summing up, the final error we have in computing $C(q,a)$ as
$\widetilde{C}(q,a)$ is
\[
  E_{final}(q,a,A,K,M,T,N)
  \leq
  \widetilde{C}(q,a)
  \frac{|E(q,a,A,M,K,N,T)|}{\Eulerphi(q)}.
\]

Practical experimentations for $q \in \{3$, \dots, $100\}$ suggested
to use different ranges for $N$ and $T$ to reach a precision of at
least $100$ decimal digits in a reasonable amount of time.
Using $A q =9600$, $M = K = 26$ and recalling that $q \mid N$ and $T$
is even, our choice is $N = (\lfloor 16800 / q \rfloor+1) q$ and
$T = 88$ if $q \in \{3$, \dots, $10\}$, while for $q \in \{90$, \dots,
$100\}$ we have to use $N = (\lfloor 40320 / q\rfloor+1) q$ and $T = 204$.
Intermediate ranges are used for the remaining integers $q$.

The programs we used to compute the Dirichlet characters $\bmod q$ and
the values of $C(q,a)$ for $q\in \{3,\dotsc,100\}$, $1\leq a\leq q$,
$(a,q)=1$, were written using the Gp scripting language of PARI/Gp;
the C program was obtained from the Gp one
using the gp2c tool. The actual computations were performed using
several LinuX pcs and one Apple MacMini computer for a total amount of
computing time equal to 1897.036096 hours = 79.043171 days.

A tiny part of the final results is collected in the following tables.
The complete set of results can be downloaded from
\url{www.math.unipd.it/~languasc/MCcomput.html} together with the
source program in Gp and the results of the verifications of the
identities \eqref{prod-over-a} and \eqref{prod-over-classes} which are
described in the section below.

\section{Verification of consistency}

The set of constants $C(q, a)$ satisfies many identities, and we
checked our results verifying that these identities hold within a very
small error.
The basic identities that we exploited are two: the first one is
\begin{equation}
\label{prod-over-a}
  \prod_{\substack{a \bmod q \\ (a, q) = 1}} C(q, a)
  =
  e^{-\gamma} \frac{q}{\Eulerphi(q)}.
\end{equation}
This can be verified using either the definition \eqref{def-C} or the
identity \eqref{value-C}, taking into account the fact that primes
dividing $q$ do not occur in any of the products $P(x; q, a)$.

The other identity is valid whenever we take two moduli $q_1$ and
$q_2$ with $q_1 \mid q_2$ and $(a, q_1) = 1$.
In this case we have
\begin{equation}
\label{prod-over-classes}
  C(q_1, a)
  =
  \prod_{\substack{j = 0 \\ (a + j q_1, q_2) = 1}}^{n - 1}
    C(q_2, a + j q_1)
  \prod_{\substack{p \mid q_2 \\ p \equiv a \bmod q_1}}
    \Bigl( 1 - \frac1p \Bigr)
\end{equation}
where $n = q_2 / q_1$.
The proof depends on the fact that the residue class $a \bmod q_1$ is
the union of the classes $a + j q_1 \bmod q_2$, for $j \in \{0$,
\dots, $ n - 1\}$.
If $q_1$ and $q_2$ have the same set of prime factors the condition
$(a + j q_1, q_2) = 1$ is automatically satisfied, since $(a, q_1) = 1$
by our hypothesis.
On the other hand, if $q_2$ has a prime factor $p$ that $q_1$ lacks,
then there are values of $j$ such that $p \mid (a + j q_1, q_2)$ and
the corresponding value of $C(q_2, a + j q_1)$ in the right hand side
of \eqref{prod-over-classes} would be undefined.
The product at the far right takes into account these primes.

To prove \eqref{prod-over-classes}, let $P(x; q, a)$ be defined by the
relation on the far left of \eqref{def-C}, without restrictions on $q$
and $a$.
Then, for $(a, q_1) = 1$ and $x \ge q_2$ write
\[
  P(x; q_1, a)
  =
  \prod_{j = 0}^{n - 1}
    P(x; q_2, a + j q_1)
  =
  \prod_{\substack{j = 0 \\ (a + j q_1, q_2) = 1}}^{n - 1}
    P(x; q_2, a + j q_1)
  \ \Pi(x; q_2, q_1, a),
\]
say.
The primes $p \le x$ such that $p \equiv a \bmod q_1$ and $p \nmid q_2$
appear in the product in the right hand side above, since there is
exactly one value of $j$ such that $p \equiv a + j q_1 \bmod q_2$ and
for any such prime it is obvious that $(a + j q_1, q_2) = 1$.
The only primes that are left are those lying in the residue class
$a \bmod q_1$ and that divide $q_2$.
Hence $\Pi(x; q_2, q_1, a)$ is exactly the product on the far right of
\eqref{prod-over-classes}.
Now \eqref{prod-over-classes} follows multiplying by a suitable
power of $\log x$ and taking the limit as $x \to +\infty$.

The validity of \eqref{prod-over-a} was checked immediately at the end
of the computation of the constants $C(q, a)$, for a fixed $q$ and for
every $1 \le a \le q$ with $(a, q) = 1$ by the same program that
computed them.
These results were collected in a file and a different program checked
that \eqref{prod-over-classes} holds within a very small error by
building every possible relation of that kind for every $q_2 \in \{3$,
\dots, $100\}$ and $q_1 \mid q_2$ with $1 < q_1 < q_2$.
The total number of identities checked is
\[
  \sum_{q = 3}^{100}
    \sum_{\substack{d \mid q \\ 1 < d < q}} \Eulerphi(d)
  =
  \sum_{q = 3}^{100} (q - 1 - \Eulerphi(q))
  =
  1907.
\]
These identities are not independent on one another, but we did not
bother to eliminate redundancies since the total time requested for
this part of the computation is absolutely negligible.
The number of independent identities is
\[
  \sum_{q = 3}^{100}
    \sum_{\substack{p \mid q \\ p < q}} \Eulerphi\Bigl( \frac qp \Bigr)
  =
  \sum_{n = 2}^{100} \pi\Bigl( \frac{100}n \Bigr) \Eulerphi(n)
  =
  1408,
\]
where $p$ denotes a prime in the sum on the left.

We wish to thank K.~Belabas and B.~Allombert for their hints about
PARI/Gp and gp2c and L.~Righi for his help using the pcs of the NumLab
laboratory \url{www.numlab.math.unipd.it} of the Department of Pure
and Applied Mathematics of the University of Padova.
We also wish to thank C.~Fassino for having read a preliminary version
of this paper.

\begin{table}
\begin{center}
\begin{tabular}{|c|c|ccc|c|}
\hline
$q$  &  $a$  &$\ $&  $C(q,a)$  &$\ $& digits \\ \hline
3  &  1  &$\ $&  1.4034774468278563951360958591826816440307\dots  &$\ $&104\\
3  &  2  &$\ $&  0.6000732161773216733074128367849176047200\dots  &$\ $&104\\
4  &  1  &$\ $&  1.2923041571286886071091383898704320653429\dots  &$\ $&104\\
4  &  3  &$\ $&  0.8689277682343238299091527791046529122939\dots  &$\ $&104\\
5  &  1  &$\ $&  1.2252384385390845800576097747492205275405\dots  &$\ $&103\\
5  &  2  &$\ $&  0.5469758454112634802383012874308140377519\dots  &$\ $&104\\
5  &  3  &$\ $&  0.8059510404482678640573768602784309320812\dots  &$\ $&104\\
5  &  4  &$\ $&  1.2993645479149779881608400149642659095025\dots  &$\ $&103\\
\vdots  & \vdots    &$\ $& \vdots  &$\ $  & \vdots    \\
9  &  1  &$\ $&  1.1738495868654491902701394683919739604995\dots  &$\ $&103\\
9  &  2  &$\ $&  0.5455303829342851960446307443914437164832\dots  &$\ $&104\\
9  &  4  &$\ $&  1.1336038613343693249917335959075962374233\dots  &$\ $&103\\
9  &  5  &$\ $&  0.9412310917798332515572574704874703583166\dots  &$\ $&103\\
9  &  7  &$\ $&  1.0547066156548587451082819988401491024340\dots  &$\ $&103\\
9  &  8  &$\ $&  1.1686623008402869661248081381642176283145\dots  &$\ $&103\\
\vdots  & \vdots    &$\ $& \vdots  &$\ $  & \vdots    \\
15 &  1  &$\ $&  1.1617073088517756555676638861655356817964\dots  &$\ $&103\\
15 &  2  &$\ $&  0.5531662836641193792434413294289420522197\dots  &$\ $&104\\
15 &  4  &$\ $&  1.1368510737193937042392719219836177668605\dots  &$\ $&103\\
15 &  7  &$\ $&  0.9888090824844727678176951687669703243697\dots  &$\ $&103\\
15 &  8  &$\ $&  1.1248826700801117041084787027689447040760\dots  &$\ $&103\\
15 & 11  &$\ $&  1.0546877248711663022320456767412694068618\dots  &$\ $&103\\
15 & 13  &$\ $&  1.0747134726382660587745323674368168616132\dots  &$\ $&103\\
15 & 14  &$\ $&  1.1429505393911402552425384830238885435764\dots  &$\ $&103\\
\vdots  & \vdots    &$\ $& \vdots  &$\ $  & \vdots    \\
21 &  1  &$\ $&  1.1141670280743936828731735756576813156065\dots  &$\ $&103\\
21 &  2  &$\ $&  0.5383301255587159174351133305605833477678\dots  &$\ $&104\\
21 &  4  &$\ $&  1.1185837991946284893102162561180399170905\dots  &$\ $&103\\
21 &  5  &$\ $&  0.8804463747350350872193530732768812838973\dots  &$\ $&103\\
21 &  8  &$\ $&  1.0809444954913878156248769107211013330026\dots  &$\ $&103\\
21 & 10  &$\ $&  1.0855302392682037293388720447231438521276\dots  &$\ $&103\\
21 & 11  &$\ $&  1.0128344672130266463968855892485398266065\dots  &$\ $&103\\
21 & 13  &$\ $&  1.0371155725767642823358797876916780548258\dots  &$\ $&103\\
21 & 16  &$\ $&  1.1035547306497255785825571380877055652196\dots  &$\ $&103\\
21 & 17  &$\ $&  1.0486412692857397440465915448981610476825\dots  &$\ $&103\\
21 & 19  &$\ $&  1.0574758123265342759648524359814135750529\dots  &$\ $&103\\
21 & 20  &$\ $&  1.1027671924237418176511972126578877947364\dots  &$\ $&103\\
\hline
\end{tabular}
\caption{\label{firsttable}
Some numerical results: the first column contains the modulus
$q$, the second the residue class $a$, the third the computed value of
$C(q,a)$ and the fourth is the number of correct decimal digits we
obtained.
The table shows the values truncated to 40 decimal digits.}
\end{center}
\end{table}

\begin{table}[htdp]
\begin{center}
\begin{tabular}{|c|c|ccc|c|}
\hline
$q$  &  $a$  &$\ $&  $C(q,a)$  &$\ $& digits \\ \hline
39 &  1  &$\ $&  1.0558043473142841979273107487867952159449\dots  &$\ $&103\\
39 &  2  &$\ $&  0.5203026628809482277529964233919621231701\dots  &$\ $&103\\
39 &  4  &$\ $&  1.0467551202397323195593324251885584436643\dots  &$\ $&103\\
39 &  5  &$\ $&  0.8477108709928609050405112584700448177533\dots  &$\ $&103\\
39 &  7  &$\ $&  0.9131634445753290856338897033232908456824\dots  &$\ $&103\\
39 &  8  &$\ $&  1.0491976120090375508070956898591030898489\dots  &$\ $&103\\
39 & 10  &$\ $&  1.0644889181790139210569905090072544013982\dots  &$\ $&103\\
39 & 11  &$\ $&  0.9611802851802015744645440449091664544815\dots  &$\ $&103\\
39 & 14  &$\ $&  1.0471282217602293552090665345631733882042\dots  &$\ $&103\\
39 & 16  &$\ $&  1.0694449785599316393966557136726680120488\dots  &$\ $&103\\
39 & 17  &$\ $&  1.0027080336857767080150127190485342860222\dots  &$\ $&103\\
39 & 19  &$\ $&  1.0063790089466405557887479935647072297591\dots  &$\ $&103\\
39 & 20  &$\ $&  1.0467993224064620442361201103591601719183\dots  &$\ $&103\\
39 & 22  &$\ $&  1.0521884311669460927257333479303503936214\dots  &$\ $&103\\
39 & 23  &$\ $&  1.0114747946261577434516887836293420101981\dots  &$\ $&103\\
39 & 25  &$\ $&  1.0597693417994788378992764465883123963780\dots  &$\ $&103\\
39 & 28  &$\ $&  1.0529671095629036217092386664444649064610\dots  &$\ $&103\\
39 & 29  &$\ $&  1.0267423753797454160121131413618162768076\dots  &$\ $&103\\
39 & 31  &$\ $&  1.0297283934645776984576326942483733223668\dots  &$\ $&103\\
39 & 32  &$\ $&  1.0482866374125031516972329668035300513497\dots  &$\ $&103\\
39 & 34  &$\ $&  1.0472581549429544593781995140831083054063\dots  &$\ $&103\\
39 & 35  &$\ $&  1.0562593819557667826211305540931669587921\dots  &$\ $&103\\
39 & 37  &$\ $&  1.0385638656749415055234884100430210797446\dots  &$\ $&103\\
39 & 38  &$\ $&  1.0674150481593719996424991312912670083485\dots  &$\ $&103\\
\hline
\end{tabular}
\caption{\label{secondtable}
Some numerical results: the first column contains the modulus
$q$, the second the residue class $a$, the third the computed value of
$C(q,a)$ and the fourth is the number of correct decimal digits we
obtained.
The table shows the values truncated to 40 decimal digits.}
\end{center}
\end{table}

\begin{table}[htdp]
\begin{center}
\begin{tabular}{|c|c|ccc|c|}
\hline
$q$  &  $a$  &$\ $&  $C(q,a)$  &$\ $& digits \\ \hline
84 &  1  &$\ $&  1.0762168747360169189445984481112147917766\dots  &$\ $&103\\
84 &  5  &$\ $&  0.8423464320992898808305526411222358430753\dots  &$\ $&103\\
84 & 11  &$\ $&  0.9670462929845278524311619985091112662169\dots  &$\ $&103\\
84 & 13  &$\ $&  0.9746953940834972813365085898448043371424\dots  &$\ $&103\\
84 & 17  &$\ $&  0.9978335235521385853486954919220491056500\dots  &$\ $&103\\
84 & 19  &$\ $&  1.0042721918535182457015722654932145385404\dots  &$\ $&103\\
84 & 23  &$\ $&  1.0128902359146896167524723309894756202191\dots  &$\ $&103\\
84 & 25  &$\ $&  1.0625109746049189658962532302336200526631\dots  &$\ $&103\\
84 & 29  &$\ $&  1.0217856732501917185719533836132834670012\dots  &$\ $&103\\
84 & 31  &$\ $&  1.0324778423499473481419749332801549343076\dots  &$\ $&103\\
84 & 37  &$\ $&  1.0448633446823406686188909998297275362347\dots  &$\ $&103\\
84 & 41  &$\ $&  1.0483511545557197512968002104563579599259\dots  &$\ $&103\\
84 & 43  &$\ $&  1.0352625518417795493214543003655548678836\dots  &$\ $&103\\
84 & 47  &$\ $&  1.0452307283367875092541042542165185077145\dots  &$\ $&103\\
84 & 53  &$\ $&  1.0473484822398583732227792995221792774100\dots  &$\ $&103\\
84 & 55  &$\ $&  1.0640407032516661060398721577715786126086\dots  &$\ $&103\\
84 & 59  &$\ $&  1.0509180585081298515408918851194537493615\dots  &$\ $&103\\
84 & 61  &$\ $&  1.0529772913206146375443030010561915545034\dots  &$\ $&103\\
84 & 65  &$\ $&  1.0629584657266981779431184953028111293016\dots  &$\ $&103\\
84 & 67  &$\ $&  1.0527738780397628191309530077617335181157\dots  &$\ $&103\\
84 & 71  &$\ $&  1.0578974865179095395282678213164071324168\dots  &$\ $&103\\
84 & 73  &$\ $&  1.0513835694502728488866738616694632680665\dots  &$\ $&103\\
84 & 79  &$\ $&  1.0561713512739106221130859861434109044623\dots  &$\ $&103\\
84 & 83  &$\ $&  1.0519063079499187595778933301342325552061\dots  &$\ $&103\\
\hline
\end{tabular}
\caption{\label{thirdtable}
Some numerical results: the first column contains the modulus
$q$, the second the residue class $a$, the third the computed value of
$C(q,a)$ and the fourth is the number of correct decimal digits we
obtained.
The table shows the values truncated to 40 decimal digits.}
\end{center}
\end{table}

\bibliographystyle{plain}

\begin{tabular}{l@{\hskip 14mm}l}
A.~Languasco               & A.~Zaccagnini \\
Universit\`a di Padova     & Universit\`a di Parma \\
Dipartimento di Matematica & Dipartimento di Matematica \\
Pura e Applicata           & Parco Area delle Scienze, 53/a \\
Via Trieste 63             & Campus Universitario \\
35121 Padova, Italy        & 43100 Parma, Italy \\
{\it e-mail}: languasco@math.unipd.it & {\it e-mail}:
alessandro.zaccagnini@unipr.it
\end{tabular}

\end{document}